\documentclass[12pt]{amsart}

\usepackage{graphicx}        
\usepackage{color}
\usepackage{amssymb,amsmath}
\usepackage{bbm}
\usepackage{ulem}  

\newtheorem{theorem}{Theorem}
\numberwithin{theorem}{section}

\theoremstyle{definition}

\newtheorem{example}[theorem]{Example}

\newcommand{\CC}{{\mathbb C}}
\newcommand{\PP}{{\mathbb P}}
\newcommand{\ZZ}{{\mathbb Z}}

\newcommand{\Fln}{{\mathbb F}\ell(n)}
\newcommand{\Fl}{{\mathbb F}\ell}

\newcommand{\Edot}{E_\bullet}
\newcommand{\Fdot}{F_\bullet}

\definecolor{navyBlue}{cmyk}{1,1,0,0.2}
\newcommand{\defcolor}[1]{{\color{blue} #1}}

\newcounter{FNC}[page]
\def\fauxfootnote#1{{\addtocounter{FNC}{2}$^\fnsymbol{FNC}$%
     \let\thefootnote\relax\footnotetext{\color{magenta}{$^\fnsymbol{FNC}$#1}}}}
\title[A geometric proof of an equivariant Pieri rule for flag manifolds]{A geometric proof of an equivariant Pieri rule\\
  for flag manifolds}
\author[C.~Li]{Changzheng Li}
\address{Changzheng Li\\
         Department of Mathematics\\
         Sun Yat-sen University\\
         China}
\email{lichangzh@mail.sysu.edu.cn}
\urladdr{http://math.sysu.edu.cn/gagp/czli}
\author[V.~Ravikumar]{Vijay Ravikumar}
\address{Vijay Ravikumar\\
         Chennai Mathematical Institute\\
         H1 SIPCOT IT Park\\
         Kelambakkam\\
         Siruseri\\
         India}
\email{vijayr@cmi.ac.in}
\urladdr{https://www.cmi.ac.in/~vijayr/}
\author[F.~Sottile]{Frank Sottile}
\address{Frank Sottile\\
         Department of Mathematics\\
         Texas A\&M University\\
         College Station\\
         Texas \ 77843\\
         USA}
\email{sottile@math.tamu.edu}
\urladdr{http://www.math.tamu.edu/~sottile}
\author[M.~Yang]{Mingzhi Yang}
\address{Mingzhi Yang\\
         Department of Mathematics\\
         Sun Yat-sen University\\
         China}
\email{yangmzh8@mail2.sysu.edu.cn}
\urladdr{}
\thanks{Research of Li and Yang supported in part by  NSFC grants  11771455 and 11831017.}
\thanks{Research of Sottile supported in part by NSF grant DMS-1501370.}

\begin{document}

\begin{abstract}
 We use geometry to give a short proof of an equivariant Pieri rule in the classical flag manifold.
 This rule is due to Robinson, who gave an algebraic proof.
\end{abstract}

\maketitle

\section*{Introduction}
 
An important problem in Schubert calculus is to find a  
formula for Schubert structure constants for flag manifolds of general Lie type.
For equivariant Schubert calculus, this is known in only two special cases, both in Lie type $A$:
the  Grassmannian, proved by Knutson and Tao \cite{KnTa}, and two-step
flag manifolds, proved by Buch \cite{Buch}.
For the manifold $\Fl(n)$ of complete flags in $\mathbb{C}^n$, a special case of this problem
is an equivariant Pieri rule that Robinson proved using algebra~\cite{Robi}.
We give a short and direct proof of this Pieri rule, using geometric arguments.

Let $\defcolor{T}\subset SL(n)$ be the diagonal torus.
The $T$-equivariant cohomology ring $H^*_T(\Fl(n))$  has an $H^*_T(\rm{pt})$-additive basis of Schubert
classes $[X_w]_T$ indexed by permutations $w$ in the symmetric group  $S_n$.
The equivariant Schubert structure constants $c_{w, v}^u$ in the product
$[X_w]_T\cdot [X_v]_T=\sum_u c_{w, v}^u[X_u]_T$ are Graham-positive~\cite{Grah}.
The {  equivariant } Pieri rule is a Graham-positive formula for $c_{w, v}^u$ when $v$ is a special permutation{, and it determines the multiplication in equivariant cohomology of any flag varety in type A.} {Its non-equivariant limit gives the classical Pieri rule   first stated by Lascoux and Sch\"utzenberger \cite{LS}.}
Our proof uses an explicit description of the projected Richardson variety associated
to $w$ and $u$ from~\cite{So}, reducing $c_{w, v}^u$ to the restriction of a special Schubert class for
a Grassmannian to a torus-fixed point as in~\cite{LiRa}.
 
The same arguments establish the same formula for torus
equivariant Chow groups of flag varieties over an algebraically closed field~\cite{Brion}.
%
%

\section{Statement of Results}

Let  $\Edot$  be the standard flag in $\CC^n$ where $E_i$ is spanned by the first $i$ standard basis
vectors, and let  $\Edot'$ be the standard opposite flag.
The Schubert variety $X_w=X_w\Edot$   of $\Fln$, where $w\in S_n$, is  defined with respect to the flag $\Edot$,
 \begin{equation}\label{Eq:SchubertVariety}
      {X_w\Edot}\ :=\
       \{ \Fdot\in\Fln \mid \dim F_i\cap E_{n+1-j}\geq\#\{k\leq i\mid w(k)\geq j\}\}\,.
 \end{equation}
This has codimension $\ell(w)$, and equals the Schubert variety of dimension $\ell(w_0w)$ with index $w_0w$
as defined in~\cite[p.\ 157]{Fu97}.
Here $\ell(w)$ is the length of the permutation $w$, and  $w_0$ is the longest permutation in $S_n$, so that
$w_0(i)=n{+}1{-}i$.
Each coefficient $c_{w, v}^u$ is either 0 or a homogeneous polynomial in
$\ZZ_{\geq 0}[t_2{-}t_1, \dotsc, t_n{-}t_{n-1}]$ of degree
$\ell(w){+}\ell(v){-}\ell(u)$ (this is Graham-positivity~\cite{Grah}).
Here, $t_i =c_1(\mathbb{C}_{\chi_i})$ is the first Chern class of the $T$-equivariant line
bundle over a point induced by  a one-dimensional representation $\mathbb{C}_{\chi_i}$ of
$T$, where $\chi_i$ denotes the character that sends $\mbox{diag}(z_1, \dotsc, z_n)\in T$   to $z_{n+1-i}$.

Fix a positive integer $m < n$.
The $m$-Bruhat order on $S_n$ is defined by $w\leq_m u$ when we have
$u=w\tau_{a_1b_1}\dotsb \tau_{a_sb_s}$ and
$\ell(w\tau_{a_1b_1}\cdots \tau_{a_ib_i})=\ell(w)+i$ for $1\leq i\leq s$.
Here,  $\tau_{a_jb_j}=(a_j,b_j)$ is a transposition with $a_j\leq m<b_j$, for $j=1,\dotsc, s$ where $s:=\ell(u)-\ell(w)$.
Consequently, when $w\leq_m u$ and $a\leq m<b$, we have $w(a)\leq u(a)$ and $w(b)\geq u(b)$.
We  write $w\xrightarrow{r_m}u$ if in addition, the integers $b_1,\dotsc, b_s$ are  distinct.

\begin{example}\label{Ex:mBO}
 When $n=9$, if $w=631594287$ in one-line notation, then
 $u=839154267=w \tau_{34}{\tau_{18}}\tau_{35}$ satisfies $w\leq_3 u$, and we also have
 $w\xrightarrow{r_3}u$ as $4,5,8$ are distinct.
 Note however that if $v=u \tau_{25}=859134267$, then $w\leq_3 v$ but we do not have $w\xrightarrow{r_3}v$.
 \hfill$\diamond$
\end{example}

Fix a positive integer $p \leq n-m$.
Let ${r(m,p)}\in S_n$ be the cyclic permutation
\[
    (m,m{+}p,m{+}p{-}1,\dotsc,m{+}1)\,.
\]
The Pieri rule involves the Schubert variety $X_{r(m,p)}$, which is defined by
a single condition,
 \[
   {X_{r(m,p)}} =\  \{\Fdot\in\Fln \mid F_m\cap E_{n+1-m-p}\neq\{0\}\}\,.
 \]
 This is the pullback of the codimension $p$ special Schubert variety in the Grassmannian $G(m,n)$ under the natural
 forgetful map.
 As the equivariant cohomology ring of any type A flag manifold is a subring of that for $\Fl(n)$, it suffices to establish
 the Pieri rule for $\Fl(n)$.

\begin{theorem}\label{Th:Main}
 For a permutation $w\in S_n$, we have the following formula in $H^*_T(\Fln)$:
 \begin{equation}
   [X_w]_T \cdot [X_{r(m,p)}]_T\ =\
    \sum_{w\xrightarrow{r_m}u}   c^u_{w,r(m,p)}  [X_u]_T\,.
 \end{equation}
 The coefficient $c^u_{w,r(m,p)}$ is nonzero if and only if $w\xrightarrow{r_m}u$ with {$r(m,p)\leq u$ and}
 $q:=p+\ell(w)-\ell(u)\geq 0$.
 When this holds,  the subset
 $\nu:=\{n{+}1{-}w(1),\dotsc,n{+}1{-}w(m)\} \cup \{n{+}1{-}w(b)\mid w(b)>u(b)\}$ consists of $m{+}p{-}q$
 elements, and defines two increasing subsequences
 \[
       {\{a_1<\dotsb< a_r\}}\ :=\
         \nu \cap \{ 1,\dotsc,n{-}m{-}p{+}1\}
 \]
  and
\[
       {\{b_1<\dotsb< b_{q+r-1}\}}\ :=\ \{n{-}m{-}p{+}2,\dotsc,n\}\smallsetminus \nu\,.
\]
 Then $c^u_{w,r(m,p)}$ is equal to $1$ if $q=0$,  or given by the following otherwise
 \begin{equation}\label{piericoeff}
    c^u_{w,r(m,p)}\ =\
    \sum\prod_{i=1}^q (t_{b_{c_i}} - t_{a_{c_i-i+1}})\ ,
 \end{equation}
 with the summation  over increasing subsequences $\{c_1<\cdots<c_q\}$ in $\{1,2,\dotsc, q{+}r{-}1\}$.
\end{theorem}
\noindent We remark that $c_{w, r(m, p)}^u$ is nonzero if and only if  the ordinary cohomology class $[X_u]$ appears in the
classical Pieri formula~\cite{LS,So} for the cup product  $[X_w]\cup [X_{r(m,p-q)}]$, which is an important part of
our proof.
The transformation $n{+}1-x$ comes from the definition~\eqref{Eq:SchubertVariety} of $X_w$.

\begin{example}\label{Ex:two}
 We compute the coefficients $c^u_{w,r(3,p)}$ for $p=3,4,5,6$, where $w$ and $u$ are from Example~\ref{Ex:mBO}.
 First, the set $\nu$ is $\{4,7,9\}\cup\{5,1,2\}$, as $w=631 594287$ and $\{5,9,8\}=\{w(b)\mid w(b)>u(b)\}$.
 (Elements of $\nu$ have the form $n{+}1{-}w(i)$.)
 When $p=3$, the coefficient $c^u_{w,r(3,3)}=1$ by the classical Pieri rule.
 When $p=4$, $q=1$ and $n{-}m{-}p{+}1=3$, so that $\{a_1<a_2\}=\{1<2\}$ and $\{b_1<b_2\}=\{6<8\}$.
 There are $2=q+r-1$ choices for $c_1$, and so we have $c^u_{w,r(3,4)}=(t_6{-}t_1) + (t_8-t_2)$.
 When $p=5$, $q=2$ and $n{-}m{-}p{+}1=2$, so that $\{a_1<a_2\}=\{1<2\}$ and $\{b_1<b_2<b_3\}=\{3<6<8\}$.
 There are three choices for $c_1<c_2$ among $\{1<2<3\}$ (namely $1<2$, $1<3$, and $2<3$), and they give
\[
   c^u_{w,r(3,5)}\ =\  (t_3{-}t_1)(t_6{-}t_1) + (t_3{-}t_1)(t_8-t_2) + (t_6{-}t_2)(t_8-t_2)\,.
\]
 Finally, when $p=6$, $q=3$ and $n{-}m{-}p{+}1=1$, so that $a_1=1$ and $\{b_1<b_2<b_3\}=\{3<6<8\}$.
 Then  $c^u_{w,r(3,6)} =  (t_3{-}t_1)(t_6{-}t_1)(t_8-t_1)$.\hfill$\diamond$
\end{example}

 \section{Proof of the equivariant Pieri rule}\label{S:two}
We prove Theorem~\ref{Th:Main} using the method of~\cite{LiRa} and exploiting the explicit description of certain
Richardson varieties and their projections in~\cite{So}.

Let \defcolor{$\Fl(1,m;n)$} denote the manifold of partial flags $F_1\subset F_m\subset\CC^n$, and
\defcolor{$G(m,n)$} denote  the  Grassmannian
of $m$-dimensional vector subspaces in $\mathbb{C}^n$.
Notice that  $\PP^{n-1}=G(1, n)$.
Let $\psi, \pi, \varphi$ denote the natural projection (forgetful) maps among these spaces.
 \begin{equation}\label{ProjMaps}
   \raisebox{-32pt}{\begin{picture}(120,70)
       \put( 90,62){$\Fln$}  \put(95,60){\vector(-3,-1){39}}
       \put(105,60){\vector(0,-1){48}} \put(107,35){$\varphi$}
               \put(25,36){$\Fl(1,m;n)$}
        \put(40,33){\vector(-1,-1){21}}  \put(50,33){\vector(2,-1){42}}
       \put(21,23){$\psi$}   \put(60,15){$\pi$}
   \put(0,0){$\PP^{n-1}$}     \put(80,0){$G(m,n)$}
  \end{picture}}
 \end{equation}
Let  \defcolor{$X^u$} denote the Schubert variety $X_{w_0u}\Edot'$,
defined with respect to the flag $\Edot'$.
Then  ${X_w^u}:=X_w\cap X^u$ is  a Richardson variety of dimension $\ell(u){-}\ell(w)$~\cite{Ric}.
Denote by $\rho^M$ the $T$-equivariant map from a $T$-space $M$ to $\rm{pt}$---a point equipped with a trivial $T$-action.
Let $Y_{(p)}$ denote the codimension $p$ special Schubert variety in $G(m,n)$.
Noticing that $X_{r(m,p)}=\varphi^{-1}Y_{(p)}=\varphi^{-1}\pi\psi^{-1}(E_{n+1-m-p})$, we compute the coefficient
$c^u_{w, r(m, p)}$ of $[X_u]_T$ in the product $[X_w]_T\cdot [X_{r(m, p)}]_T$ using the equivariant push forward to a
point and the projection formula,
  \begin{align*}
      c^u_{w, r(m, p)}&= \rho_*^{\Fl(n)}\big( [X_w]_T\cdot [X_{r(m, p)}]_T\cdot[X^u]_T\big)\
                        =\ \rho_*^{\Fl(n)}\big([X_w^u]_T\cdot [X_{r(m, p)}]_T\big) \\
               &=\rho_*^{G(m,n)}\big(\varphi_*[X_w^u]_T\cdot [Y_{(p)}]_T\big)\\
               &=\rho_*^{\Fl(1,m;n)}\big(\pi^*\varphi_*[X_w^u]_T\cdot \psi^*[E_{n+1-m-p}]_T\big)\\
               &=\rho_*^{\PP^{n-1}} \left( ( \psi_* \pi^* \varphi_* [X_w^u]_T)\cdot[E_{n+1-m-p}]_T\right)\ .
  \end{align*}

Let ${Y^u_w}:=\varphi(X^u_w)$ be the image of $X^u_w$ in $G(m,n)$ and
${Z^u_w}:=\psi\circ\pi^{-1}(Y^u_w)$, which is the set of points in $\PP^{n-1}$ that lie on some $m$-plane in a
flag in $X^u_w$.
%
%
We have that $\dim Y^u_w\leq \dim X^u_w$ and $\dim Z^u_w\leq m + \dim Y^u_w$, since
the general fiber of $\pi$ has dimension $m$.
Moreover, the class
$\psi_* \pi^* \varphi_* [X_w^u]_T$ is zero if either inequality is strict.

The Chevalley formula~\cite{FuLa} expresses  $[X_w]_T\cdot [X_{r(m, 1)}]_T$ as a Graham-positive sum of classes
$[X_u]_T$, where either $u=w$ or $u$ covers $w$ in the $m$-Bruhat order.
Iterating the Chevalley formula shows that $[X_{r(m, p)}]_T$ is a term of $([X_{r(m, 1)}]_T)^p$.
Thus $[X_w]_T\cdot [X_{r(m, p)}]_T$  is a subsum of $[X_w]_T\cdot ([X_{r(m, 1)}]_T)^p$, which is a Graham-positive
combination of classes $[X_u]_T$ with $w\leq_m u$, again by the Chevalley formula.
Thus $c^u_{w, r(m, p)}=0$ unless $w\leq_m u$.
Assuming this, $\varphi_* [X_w^u]_T$ is a positive multiple of $[Y_w^u]_T$.
By~\cite[Lemma 10]{So}, $Z^u_w$ is a
subset of a linear subspace (witten there as $Y$) of $\PP^{n-1}$ of dimension $m+\#\{m<b\mid w(b)> u(b)\}$.
From the definition of the $m$-Bruhat order and $\xrightarrow{r_m}$, this   is strictly less than
$m+\dim X_w^u$, unless $w\xrightarrow{r_m} u$.
Thus $c^u_{w,r(m,p)}\neq 0$ only if $w\xrightarrow{r_m} u$.

We recall Lemma~15 of~\cite{So}, which identifies both $Y_w^u$ and $Z^u_w$ when $w\xrightarrow{r_m}u$
(and shows that the maps $\varphi$ and $\psi$ to them are birational).
This differs from the statement in~\cite{So} in that
our standard basis is different from the basis used there{, with $e_i$ here being $e_{n+1-i}$ in~\cite{So}}.\medskip

\noindent{\bf Lemma 15} from~\cite{So}.\
 {\it Suppose that $w\xrightarrow{r_m} u$ and $u=w\tau_{a_1 b_1}\dotsb\tau_{a_s b_s}$ with
 $a_i\leq m < b_i$ and $\ell(w\tau_{a_1 b_1}\dotsb \tau_{a_i b_i})=\ell(w)+i$ for $1\leq i\leq s=\ell(u)-\ell(w)$.
 Define
 \begin{eqnarray*}
   L_j &=& \langle e_{n+1-w(j)}, e_{n+1-w(b_i)}\mid a_i=j\rangle\qquad j=1,\dotsc,m\\
   M   &=& \langle e_{n+1-w(k)} \mid m<k \mbox{ and } w(k)=u(k)\rangle\,.
 \end{eqnarray*}
 Then $\dim L_j=1+\#\{i\mid a_i=j\}$, if $\Fdot\in X^u_w$ then $\dim F_m\cap L_j=1$ for
 $1\leq j\leq k$.
 Consequently, the image $Y^u_w$ is a Richardson variety in the Grassmannian with respect to different coordinate flags than
 $\Edot$ and $\Edot'$, and the map $\varphi\colon X^u_w\to Y^u_w$ has degree $1$.}\medskip

(The $e_{n+1-w(j)}$ in the definition of $L_j$ corrects a typographic error in the published article.)
As the map $\varphi\colon X^u_w\to Y^u_w$ has degree $1$, we have $\varphi_*[X^u_w]_T=[Y^u_w]_T$.
More can be said about $Y^u_w$, but the important consequence for now is that
$[Z^u_w]_T= \psi_* \pi^* \varphi_* [X_w^u]_T$, and
\[
  Z^u_w\ =\ L_1\oplus\dotsb\oplus L_m\ =\
  \langle e_{n+1-w(j)}, e_{n+1-w(b_i)}\mid j=1,\dotsc,m\,,\ i=1,\dotsc,r\rangle\ =\  \langle e_k \mid k \in \nu \rangle\,.
\]
 Here, $\nu$ is the indexing set defined in the statement of Theorem \ref{Th:Main}.
 As in \cite{LiRa}, $Z^u_w$ equals the projected Richardson variety coming from the single point Richardson variety
 $\langle e_k \mid k \in \nu \rangle$ in the Grassmannian $G(m{+}p{-}q,n)$, via the standard maps~\eqref{ProjMaps}
 from $\Fl(1,m{+}p{-}q;n)$.
 By~\cite[Proposition 4.2]{LiRa}, $c_{w, r(m, p)}^u$ equals $N^\nu_{\nu,q}$, the localization of the special Schubert class
 of codimension $q$ in the Grassmannian $G(m{+}p{-}q,n)$ at the torus-fixed point $\langle e_k \mid k \in \nu \rangle$.
 The indexing set $\nu$   is   a Schubert symbol in \cite{LiRa}, and the Grassmannian permutation associated to
 $\nu$ is obtained by sorting  $\tilde \nu:=\{w(1),\ldots, w(m)\}\cup \{w(b)\mid w(b)>u(b)\}$
 in increasing order.
 Thus $N_{\nu, q}^{\nu}$ is nonzero  if and only if   $q+ m+p-q\leq \max \tilde \nu=\max \{u(1),\ldots, u(m)\}$, namely
 $r(m,p)\leq u$ in the Bruhat order for $S_n$. 
 Here   the last equality holds due to the observations: (i) $u(a_i)=w(b_s)>u(b_s)$ with $s=\max\{k\mid a_i=a_k, k\leq r\}$
 so that $u(j)\in \tilde \nu$ for any $1\leq j\leq m$;  (ii)   $w(j)\leq u(j)$ for $1\leq j\leq m$, and if $w(b)>u(b)$ then
 $b=b_i$ for a unique $1\leq i\leq m$, implying that $u(a_i)\geq w(b)$ since $w\xrightarrow{r_m}u$. 
By~\cite[Appendix A]{LiRa}   {the nonvanishing} localization is~\eqref{piericoeff}, which completes the proof.

From the proof of Lemma~15 in~\cite{So}, there are partitions $\mu\subset\lambda$ for $G(m,n)$, and a permutation
$\omega\in S_n$ such that $c^u_{w,r(m,p)}=\omega(c^\lambda_{\mu,(p)})$.

\begin{example}
We illustrate Lemma~15.
An invertible matrix $F$ gives a flag $\Fdot$, where $F_i$ is the row space of the first $i$ rows of $F$.
The two matrices below represent general elements of $X_w$ and $X^u$, for $w$ and $u$ from Example~\ref{Ex:mBO}.
Here, $\cdot$ is zero and $*\in\CC$ is arbitrary. 
\[
   \left(\begin{matrix}
     *  &  *  &  *  &  1  &\cdot&\cdot&\cdot&\cdot&\cdot\\
     *  &  *  &  *  &\cdot&  *  &  *  &  1  &\cdot&\cdot\\
     *  &  *  &  *  &\cdot&  *  &  *  &\cdot&  *  &  1  \\
     *  &  *  &  *  &\cdot&  1  &\cdot&\cdot&\cdot&\cdot\\
     1  &\cdot&\cdot&\cdot&\cdot&\cdot&\cdot&\cdot&\cdot\\
   \cdot&  *  &  *  &\cdot&\cdot&  1  &\cdot&\cdot&\cdot\\
   \cdot&  *  &  *  &\cdot&\cdot&\cdot&\cdot&  1  &\cdot\\
   \cdot&  1  &\cdot&\cdot&\cdot&\cdot&\cdot&\cdot&\cdot\\
   \cdot&\cdot&  1  &\cdot&\cdot&\cdot&\cdot&\cdot&\cdot
  \end{matrix}\right)
   \qquad
  \left(\begin{matrix}
    \cdot&  1  &  *  &  *  &  *  &  *  &  *  &  *  &  *  \\
    \cdot&\cdot&\cdot&\cdot&\cdot&\cdot&  1  &  *  &  *  \\
      1  &\cdot&  *  &  *  &  *  &  *  &\cdot&  *  &  *  \\
    \cdot&\cdot&\cdot&\cdot&\cdot&\cdot&\cdot&\cdot&  1  \\
    \cdot&\cdot&\cdot&\cdot&  1  &  *  &\cdot&  *  &\cdot\\
    \cdot&\cdot&\cdot&\cdot&\cdot&  1  &\cdot&  *  &\cdot\\
    \cdot&\cdot&\cdot&\cdot&\cdot&\cdot&\cdot&  1  &\cdot\\
    \cdot&\cdot&\cdot&  1  &\cdot&\cdot&\cdot&\cdot&\cdot\\
    \cdot&\cdot&  1  &\cdot&\cdot&\cdot&\cdot&\cdot&\cdot\\
  \end{matrix}\right)
\]
The Richardson variety $X^u_w$ has a parametrization by triples $(\alpha,\beta,\gamma)$ of nonzero complex numbers.
We show this in two equivalent ways.
The matrix on the left lies in $X_w$ and that on the right in $X^u$.
For every $i=1,\dotsc,9$ the first $i$ rows of both matrices have the same span.
\[
  \left(\begin{matrix}
    \cdot&\beta&\cdot&  1  &\cdot &\cdot&\cdot&\cdot&\cdot\\
    \cdot&\cdot&\cdot&\cdot&\cdot &\cdot&  1  &\cdot&\cdot\\
   \gamma&\cdot&\cdot&\cdot&\alpha&\cdot&\cdot&\cdot&  1  \\\hline
   \gamma&\cdot&\cdot&\cdot&\alpha&\cdot&\cdot&\cdot&\cdot\\
   \gamma&\cdot&\cdot&\cdot&\cdot &\cdot&\cdot&\cdot&\cdot\\
    \cdot&\cdot&\cdot&\cdot&\cdot &  1  &\cdot&\cdot&\cdot\\
    \cdot&\cdot&\cdot&\cdot&\cdot &\cdot&\cdot&  1  &\cdot\\
    \cdot&\beta&\cdot&\cdot&\cdot &\cdot&\cdot&\cdot&\cdot\\
    \cdot&\cdot&  1  &\cdot&\cdot &\cdot&\cdot&\cdot&\cdot
  \end{matrix}\right)
  \qquad
  \left(\begin{matrix}
    \cdot&\beta&\cdot&  1  &\cdot &\cdot&\cdot&\cdot&\cdot\\
    \cdot&\cdot&\cdot&\cdot&\cdot &\cdot&  1  &\cdot&\cdot\\
   \gamma&\cdot&\cdot&\cdot&\alpha&\cdot&\cdot&\cdot&  1  \\\hline
    \cdot&\cdot&\cdot&\cdot&\cdot &\cdot&\cdot&\cdot&  1  \\
    \cdot&\cdot&\cdot&\cdot&\alpha&\cdot&\cdot&\cdot&  1  \\
    \cdot&\cdot&\cdot&\cdot&\cdot &  1  &\cdot&\cdot&\cdot\\
    \cdot&\cdot&\cdot&\cdot&\cdot &\cdot&\cdot&  1  &\cdot\\
    \cdot&\cdot&\cdot&  1  &\cdot &\cdot&\cdot&\cdot&\cdot\\
    \cdot&\cdot&  1  &\cdot&\cdot &\cdot&\cdot&\cdot&\cdot
  \end{matrix}\right)
\]
The row span of their first three rows parameterizes the projected Richardson variety $Y^u_w$.
Since the remaining rows depend upon the first three, the map $X^u_w\to Y^u_w$ is birational.\hfill$\diamond$
\end{example}

\providecommand{\bysame}{\leavevmode\hbox to3em{\hrulefill}\thinspace}
\providecommand{\MR}{\relax\ifhmode\unskip\space\fi MR }
\providecommand{\MRhref}[2]{%
  \href{http://www.ams.org/mathscinet-getitem?mr=#1}{#2}
}
\providecommand{\href}[2]{#2}

\end{document}